\documentclass[a4paper,12pt]{article}

\usepackage{graphicx, amssymb, latexsym, amsfonts, amsmath, lscape, amscd,
	amsthm, color, epsfig}

\usepackage[T1]{fontenc}
\usepackage[latin1]{inputenc}

\usepackage[babel=true]{csquotes}

\usepackage{authblk}
\usepackage{setspace}

\usepackage{enumitem}
\usepackage{tikz-network}
\usepackage{tikz}
\usepackage{mathtools}

\DeclarePairedDelimiter\floor{\lfloor}{\rfloor}
\usetikzlibrary{decorations.pathreplacing}
\newtheorem{theorem}{Theorem}
\numberwithin{theorem}{section}
\newtheorem{conjecture}{Conjecture}
\numberwithin{conjecture}{section}
\newtheorem{corollary}{Corollary}
\numberwithin{corollary}{section}

\numberwithin{example}{section}
\newtheorem{lemma}{Lemma}
\numberwithin{lemma}{section}

\numberwithin{proposition}{section}

\numberwithin{problem}{section}

\numberwithin{remark}{section}
\newtheorem{claim}{Claim}

\numberwithin{definition}{section}
\numberwithin{equation}{section}
\date{}

\begin{document}
	
	\renewcommand\qedsymbol{$\blacksquare$}
	
	\newcommand\subrel[2]{\mathrel{\mathop{#2}\limits_{#1}}}
	
	\author {Abbas Alhakim$^2$ Mouhamad El Joubbeh $^{1,2}$}
	
		\footnotetext[1]{American University of Beirut, Department of mathematics , Beirut, Lebanon.}
	\footnotetext[2]{Lebanese University, KALMA Laboratory, Baalbeck, Lebanon.}

	\title{Subdivisions of Oriented Cycles in Digraphs with Hamiltonian directed path}
	
	\maketitle

	\begin{abstract}
		Cohen et al. conjectured that for every oriented cycle $C$ there exist an integer $f(C)$ such that every strong $f(C)$-chromatic digraph contains a subdivision of $C$. El Joubbeh confirmed this conjecture for Hamiltonian digraphs. Indeed, he showed that every $3n$-chromatic Hamiltonian digraph contains a subdivision of every oriented cycle of order $n$. In this article, we improve this bound to $2n$. Furthermore, we show that, if $D$ is a digraph containing a Hamiltonian directed path with chromatic number at least $12n-5$, then $D$ contains a subdivision of every oriented cycle of order $n$. Note that a digraph containing a Hamiltonian directed path need not be strongly connected. Thus, our current result provides a deeper understanding of the condition that may be needed to fully solve the conjecture. 
	\end{abstract}
	\section{Introduction}
	\noindent 
	In graph theory, we represent a graph as $G = (V(G), E(G))$, where $V(G)$ signifies the set of vertices (nodes) in the graph, and $E(G)$ denotes the set of edges, representing the connections between the vertices. In this paper, we assume that all graphs are simple, meaning they lack loops and multiple edges. A \textit{vertex enumeration} of a graph is a listing or ordering of its vertices. It is represented as a sequence or set of vertices, typically denoted as $N=(v_{i_1}, \dots, v_{i_n})$, where $V(G)=\{v_1,\dots,v_n\}$ and $\{i_1,\dots,i_n\}=\{1,\dots,n\}$. A \textit{Hamiltonian path} (or \textit{Hamiltonian cycle}) of a graph refers to a path (or cycle) that traverses all its vertices. The \textit{length} of a path $P$ (or cycle $C$) is represented as $l(P)$ (or $l(C)$) and indicates the number of edges in that path or cycle. For a given graph $G$, the \textit{girth} of $G$, denoted as $g(G)$, corresponds to the length of the shortest cycle within $G$. A \textit{proper $k$-coloring} of a graph $G$ is defined as a function $c:V(G)\to\{1,\dots,k\}$ such that adjacent vertices possess distinct colors. A graph is called \textit{$k$-colorable} if it admits a proper $k$-coloring. The \textit{chromatic number} of a graph, denoted as $\chi(G)$, represents the smallest integer $k$ for which $G$ is $k$-colorable. A graph with a chromatic number of $k$ is referred to as a \textit{$k$-chromatic graph}. A \textit{spanning subgraph} of $G$ is a subgraph $H$ of $G$ such that $V(H)=V(G)$. An \textit{induced subgraph} $H$ of a graph $G$ is defined as a subgraph that includes all edges of $G$ connecting two vertices within the set $V(H)$. Formally, we write $H=G[V(H)]$.\\	
	
	An \textit{oriented graph}, or simply a \textit{digraph}, is formed by assigning directions to the edges of a given graph. The notation $D = (V(D), E(D))$ is a concise way of describing the structure of a digraph $D$, where $V(D)$ is the set of vertices of $D$ and $E(D)$ is the set of directed edges (\textit{arcs}) of $D$, where each arc is an ordered pair of vertices $(u, v)$ indicating a directed connection from vertex $u$ to vertex $v$. Conversely, the graph derived from a digraph $D$ by ignoring its edge directions is referred to as the \textit{underlying graph} of $D$, denoted as $G(D)$. As such, the principles established for graphs can be seamlessly extended to digraphs by directly applying them to the underlying graphs of those digraphs. A \textit{directed path}, is a finite sequence of distinct vertices and directed edges, where each edge connects one vertex to the next in the sequence, and all edges have a consistent direction. A \textit{directed cycle}, or simply a \textit{circuit}, is a closed directed path, forming a loop by returning to the starting vertex. A \textit{Hamiltonian digraph} is one that contains a Hamiltonian circuit. A \textit{block} of a non-directed oriented cycle $C$ (respectively, path $P$) is a maximal directed sub-path of $C$ (respectively, of $P$). Assume that a non-directed oriented cycle $C$ is composed of $l$ blocks of consecutive lengths $k_{1}, \dots , k_{l}$, then we denote this cycle as $C=C^+(k_{1}, \dots , k_{l})$ if the block of length $k_1$ is forward, and as $C=C^-(k_{1}, \dots , k_{l})$ otherwise. It is evident, by the definition of a block, that $l$ must be an even number. It is worth noting that, in general, $C^+(k_1,\dots,k_l) \neq C^-(k_1,\dots,k_l)$; for instance, $C^+(1,1,2,2)\neq C^-(1,1,2,2)$. A non-directed oriented cycle $C$ is considered \textit{antidirected} if all of its blocks are of length $1$. A \textit{subdivision} of a digraph $D$ is a digraph obtained from $D$ by replacing each arc oriented from $x$ to $y$ with an $x,y$-directed path. In a circuit $C$, if we consider two vertices $x$ and $y$, then $C_{[x,y]}$ denotes the unique directed sub-path of $C$ beginning at $x$ and ending at $y$. Similarly, in a directed path $P$ with vertices $x$ and $y$, where there exists a subpath within $P$ starting from $x$ and ending at $y$, the notation $P_{[x,y]}$ represents the directed sub-path of $P$ from $x$ to $y$. If $(x, y)$ is an arc in a directed path $P$ or a circuit $C$, $y$ is referred to as the \textit{successor} of $x$ on $P$ or in $C$, respectively.\\
	
An interesting problem in graph theory revolves around identifying an integer $k$ such that every $k$-chromatic digraph includes a copy of each member of a specified family of digraphs. For instance, Burr \cite{Burr01} established that a digraph with chromatic number $(n-1)^2$ includes every oriented tree of order $n$ and conjectured that this bound might be refined to $2n-2$.\ While having a large chromatic number may stem from including a substantial complete subgraph, Erd\"os demonstrated in \cite{girth} that for any positive integers $g$ and $k$, there exists a graph $G$ with girth exceeding $g$ and chromatic number exceeding $k$.\ This implies that if an element within a given family of digraphs $\cal D$ has a bounded girth, finding an integer $k$ for which every $k$-chromatic digraph includes a representation of each element in $\cal D$ becomes impossible. This shifts the focus of the problem to the study of two distinct classes of digraph families; those whose elements contain no oriented cycles (trees and paths), and those whose elements have no upper bound on their girths (subdivision of cycles and bispindles).\\

Concerning the oriented trees, Addario-Berry et al.\ \cite{Addario-Berry01} have improved Burr's bound by showing that every oriented tree of order $n$ is contained in every $(\frac{n^2}{2}-\frac{n}{2}+1)$-chromatic digraph, which is the best bound reached so far.\ Regarding paths, the most celebrated result, known as Gallai-Hasse-Roy-Vitaver theorem  \cite{Gallai01,Hasse,Roy01,Vitaver}, deals with directed paths and states that every $n$-chromatic digraph contains a directed path of length $n-1$.\ El Sahili \cite{sahili2} conjectured that every oriented path of order $n$ exists in every $n$-chromatic digraph.\ Addario et al.\ \cite{Addario-Berry01} and El Sahili et al.\ \cite{sahili1}, independently proved this conjecture for two-block paths with $n\geq 4$.\ For paths with $t$ blocks, $t\geq 3$, El joubbeh et al. \cite{El.joubbeh,El.joubbeh.Ghazal} showed that every oriented $n$-path with $t$ blocks is contained in every $(4^{r(t)}n+q(t))$-chromatic digraph where $r(t)=\floor*{\log_2(t-1)}$ and  $q(t)=-\frac{20}{7}8^{r(t)-1}+\frac{6}{7}$.\ Additional results were presented for three-block paths with certain length restrictions on their blocks (See \cite{Mortada2, Mortada1, Tarhini.2, Tarhini}).\\
	
On the flip side, Bondy \cite{bondy} showed that each strong digraph has a circuit that is longer than its chromatic number. It's crucial to note that the assumption of strong connectivity is indispensable because an acyclic $k$-chromatic digraph lacks any circuit. Interestingly, even when contemplating the presence of non-directed oriented cycles, the generality cannot be extended to general digraphs, due to the existence of acyclic digraphs with arbitrarily large chromatic numbers and no oriented cycles of two blocks, as observed by Gy\'{a}rf\'{a}s and Thomassen (see \cite{Addario-Berry01}). Cohen et al. \cite{cohen1} expanded this observation to any number of blocks. Specifically, they demonstrated that for any positive integers $b$ and $c$, there exists an acyclic digraph $D$ with $\chi(D) \geq c$, where all oriented cycles in $D$ have more than $b$ blocks. This underscores the importance of confining the exploration of oriented cycles in chromatic digraphs to a specific class of digraphs. Consequently, Cohen et al.\ \cite{cohen1} suggested that a ``strong condition'' must be imposed on digraphs to enhance the possibility of ensuring the existence of a subdivision of the desired cycle.
\begin{conjecture}\label{conj}
For any oriented cycle $C$, there exists a constant $f(C)$ such that every $f(C)$-chromatic strong digraph contains a subdivision of $C$.
\end{conjecture} 
In their work, Cohen et al. \cite{cohen1} substantiated Conjecture \ref{conj} for the case of two-block cycles. In fact, they demonstrated that the chromatic number of strong digraphs with no subdivisions of two-block cycles $C^+(k_{1},k_{2})$ is bounded above by $(k_{1}+k_{2}-2)(k_{1}+k_{2}-3)(2k_{2}+2)(k_{1}+k_{2}+1)$. The most recent upper bound, as determined by Kim et al. \cite{Kim}, is $12k^2$, where $k=\max\{k_1,k_2\}$. For cycles with four blocks, Al-Mniny \cite{k111} showed that every $D$ is a strong digraph with no subdivision of $C^+(k,1,1,1)$ implies $\chi(D)\leq 8^3k$. Mohsen \cite{k111mohsen} improved this bound to $18k$. Al-Mniny et al. \cite{mnza} demonstrated that if $D$ is a strong digraph with no subdivision of $C^+(k,1,k,1)$, then $\chi(D)\leq 36(2k)(4k+2)$. El Joubbeh \cite{El.joubbeh.cycles} later confirmed Conjecture \ref{conj} for Hamiltonian digraphs, showing that every $3n$-chromatic Hamiltonian digraph contains a subdivision of every oriented cycle of order $n$.\\
 
 A similar problem concerns the $(2+1)$-bispindle. Such a bispindle, denoted as $B(k_1,k_2;k_3)$, consists of two directed paths from $x$ to $y$ with lengths $k_1$ and $k_2$ respectively, along with one directed path from $y$ to $x$ of length $k_3$, and all these directed paths being pairwise internally disjoint. Cohen et al. \cite{bis-cohen} conjectured that for any bispindle $B(k_1,k_2;k_3)$, there is an integer $g(k_1,k_2,k_3)$ such that every $g(k_1,k_2;k_3)$-chromatic strong digraph contains a subdivisions of $B(k_1,k_2;k_3)$. In the same paper, they proved this conjecture for specific cases like $B(k_1,1;k_3)$, $B(k_1,1;1)$, and $B(1,k_2;1)$. Additionally, Al-Mniny \cite{bis-mniny} confirmed this conjecture for $B(k_1,k_2;1)$.\\

Surprisingly, in this paper, we obtained a result aligned with Conjecture \ref{conj} but without the ``strong condition''. Instead, we imposed another condition on digraphs, stating that if a digraph contains a Hamiltonian directed path and has a chromatic number of $12n-5$, then it contains a subdivision of every oriented cycle of order $n$. Notably, such a digraph can be non-strong or even acyclic, providing a more profound insight into the necessary conditions for fully resolving the conjecture. Further findings in this paper are presented in Sections \ref{secant} and \ref{H.C.}, which are necessary before unveiling our main result in Section \ref{H.P.}.
	\section{$k$-secant edges} \label{secant}
	Consider a graph $G$ equipped with a vertex enumeration $N=(v_1, \dots, v_t)$. Al Mniny et al.\ \cite{k111} introduced the concept of \textit{secant edges} with respect to $N$, defined as two edges $v_iv_j$ and $v_rv_l$ in $G$, where $i<r<j<l$. El Joubbeh \cite{El.joubbeh.cycles} extended this notion to $k$-\textit{secant edges}, defining pairs $v_iv_j$ and $v_rv_l$ as $k$-\textit{secant edges} with respect to $N$ if they are secant and $j-r\geq k$.\ El Joubbeh \cite{El.joubbeh.cycles} showed that if $G$ contains no $k$-secant edges w.r.t. $N$, then $G$ is $(4k-1)$-colorable.\ Subsequently, Al Mniny et al.\ \cite{secant} showed that it is $3k$-colorable. In this paper, we improve this bound to $2k+1$, and put forth a conjecture suggesting its $(k+2)$-colorability.
	\begin{lemma}\label{L:secant}
	Let $G$ be a graph with an enumeration $N$. If $G$ contains no $k$-secant edges with respect to $N$, then $G$ is $(2k+1)$-colorable.
\end{lemma}

\begin{proof}
	We use induction on $v(G)$. If $v(G)\leq 2k+1$, then it is trivial. Now, let $G$ be a graph with an enumeration $N=(v_1,\ldots,v_s)$ that contains no $k$-secant edges with respect to $N$. 
	\begin{claim}\label{claim1}
		$G$ has a vertex $v_{i_0}$ such that $d_G(v_{i_0})\leq 2k$.
	\end{claim} 

		Proof of Claim 1. If not, then $d_G(v)\geq 2k+1$, $\forall v\in V(G)$. Thus, for all $v\in V(G)$, we can choose a neighbor $\tilde{v}$ of $v$ such that $|i_N(v)-i_N(\tilde{v})|\geq k+1$, where $i_N(v)$ is the index of $v$ in $N$.\ We define the following sets:\\
		$$right(G)=\{v\in V(G) \ | \ i_N(\tilde{v})>i_N(v)\}$$
		$$left(G)=\{v\in V(G) \ | \ i_N(\tilde{v})<i_N(v)\}.$$
		Note that $right(G)$ and $left(G)$ contain the first $k$ vertices and the last $k$ vertices of $N$, respectively.  Thus the two sets form a partition of $V(G)$. Let $v_j\in right(G)$ such that $j$ is maximal. It follows by this maximality that $v_{j+k}\in left(G)$. Moreover, by the defining condition of $\tilde{v}$, the index of $\tilde{v}_{j+k}$ in $N$ is less than $j$. Therefore,
		the two edges $v_j\tilde{v}_j$ and $v_{j+k}\tilde{v}_{j+k}$ are $k$-secant edges (see Figure \ref{fig:path02}). This contradiction proves the claim \ref{claim1}.
		\begin{figure}[h] 
			\centering
			\caption{}
			
			\begin{tikzpicture}
				\Vertex[x=1.5,size=0.2,color=black,label={$v_1$},position=below]{v1}
				\Vertex[x=5,size=0.2,color=black,label={$v_i$},position=below]{v3}
				\Vertex[x=8.8,size=0.2,color=black,label={${\tilde{v}}_i$},position=below]{v4}
				\Vertex[x=10.5,size=0.2,color=black,label={$v_s$},position=below]{v5}
				\Vertex[x=7.5,size=0.2,color=black,label={$v_{i+k}$},position=above]{v6}
				\Vertex[x=3.5,size=0.2,color=black,label={${\tilde{v}}_{i+k}$},position=above]{v7}
				\Edge[style=dashed,color=brown](v1)(v3)
				\Edge[style=dashed,color=brown](v6)(v4)
				\Edge[color=brown](v3)(v6)
				\Edge[style=dashed,color=brown](v4)(v5)
					\Edge[bend=35,color=black](v3)(v4)
					\Edge[bend=30,color=black](v6)(v7)
				
				\draw[color=brown] (2.25,0.2) to[out=45, in=135] (2.75,0.2);
				\draw[color=brown] (9.25,0.2) to[out=45, in=135] (9.75,0.2);
				\foreach \x in {2.25,2.5,2.75,5.67, 6,6.33,6.67,7,9.25,9.5,9.75}{
					\filldraw[black] (\x,0) circle (2pt);
				}
			
			\end{tikzpicture}
			
			\label{fig:path02}
		\end{figure}\\

		Let $G'=G-\{v_{i_0}\}$ and $N'=(v_1,\ldots,v_{i_0-1},v_{i_0+1},\ldots ,v_s)$. Note that every pair of $k$-secant edges in $G'$ with respect to $N'$ is also a pair of $k$-secant edges in $G$ with respect to $N$. This forces $G'$ to have no $k$-secant edges either. Provoking the induction hypothesis, we get  a proper $2k+1$-coloring $c'$ of $G'$. Since $d(v_{i_0})\leq2k$, $c'$ can be extended to a proper $2k+1$-coloring $c$ of $G$.

\end{proof}
A direct corollary follows from Lemma \ref{L:secant}.
\begin{corollary}\label{coro1}
	Let $G$ be a digraph with a vertex enumeration $N$, and let $k$ be a positive integer. If $\chi(G) \geq 2k+2$, then $G$ contains $k$-secant edges with respect to $N$.
\end{corollary}

Consider a digraph $D$ with a vertex enumeration $N$. Two arcs $e$ and $f$ in $D$ are $k$-\textit{secant arcs} with respect to $N$ if their corresponding edges in $G(D)$ are $k$-secant edges with respect to $N$. Thus, the two previous results apply directly to digraphs.\\ 

\section{Subdivision of oriented cycles in Hamiltonian digraphs}\label{H.C.}
		This section focuses on refining the chromatic number bound for Hamiltonian digraphs that include a subdivision of every oriented cycle of order $n$, aiming to lower the bound from $3n$ to $2n$. In his article, El Joubbeh \cite{El.joubbeh.cycles} showed that if $C=C^+(k_1,\dots, k_{2t})$ and $D$ a Hamiltonian digraph with $\chi(D)\geq (k_1+1)+4k_2\dots+ (k_{2t-1}+1)+4k_{2t}$ then $D$ contains a subdivision of $C$. His proof involves constructing a list of induced sub-digraphs $(D_i)_{1\leq i\leq 2t}$ of $D$ with specific properties. Starting with a Hamiltonian circuit $C'$ of $D$, each sub-digraph $D_i$ is created such that its chromatic number alternates between $k_i+1$ when $i$ is odd and $4k_i$ when $i$ is even, and $C'[V(D_i)]$ forms a Hamiltonian directed path in $D_i$ in a way that the start vertex of $C'[V(D_i)]$ for $i > 1$ is the successor of the end vertex of $C'[V(D_{i-1})]$ (See Figure \ref{fig:path0}).
		\begin{figure}[h]
			\centering
			\caption{}
	\begin{tikzpicture}
		\draw[line width=1pt, brown] (0,0) circle (2);
			\node[align=center,rotate=-90, brown] at (0,2) {\footnotesize{$\blacktriangle$}};
		\pgfmathsetmacro{\xcoord}{2*cos(30)}
		\pgfmathsetmacro{\ycoord}{2*sin(30)}
		\node[align=center, rotate=-150, brown] at (\xcoord,\ycoord) {\footnotesize{$\blacktriangle$}};
		\pgfmathsetmacro{\xcoord}{2*cos(150)}
		\pgfmathsetmacro{\ycoord}{2*sin(150)}
		\node[align=center, rotate=-30, brown] at (\xcoord,\ycoord) {\footnotesize{$\blacktriangle$}};
		\pgfmathsetmacro{\xcoord}{2*cos(-33)}
		\pgfmathsetmacro{\ycoord}{2*sin(-33)}
		\node[align=center, rotate=-213, brown] at (\xcoord,\ycoord) {\footnotesize{$\blacktriangle$}};
		\pgfmathsetmacro{\xcoord}{2*cos(-80)}
		\pgfmathsetmacro{\ycoord}{2*sin(-80)}
		\node[align=center, rotate=-260, brown] at (\xcoord,\ycoord) {\footnotesize{$\blacktriangle$}};
		\pgfmathsetmacro{\xcoord}{2*cos(-120)}
		\pgfmathsetmacro{\ycoord}{2*sin(-120)}
		\node[align=center, rotate=-300, brown] at (\xcoord,\ycoord) {\footnotesize{$\blacktriangle$}};
		\node[align=center] at (2.3,1.2) {$D_1$};
		\node[align=center] at (-2.3,1.2) {$D_{2t}$};
		\node[align=center] at (2,-1.65) {$D_2$};
		\node[align=center] at (-1.5,-2) {$D_i$};
		\pgfmathsetmacro{\xcoord}{2*cos(-160)}
		\pgfmathsetmacro{\ycoord}{2*sin(-160)}
		\node[align=center, rotate=-340, brown] at (\xcoord,\ycoord) {\footnotesize{$\blacktriangle$}};
		\foreach \i in {1,2,3,4} {
			\draw ({360/6 * (\i - 1)}:2) node[fill, circle, inner sep=1.75pt,] {};
		}
		\draw (-7:2) node[fill, circle, inner sep=1.75pt,] {};
		\draw (-60:2) node[fill, circle, inner sep=1.75pt,] {};
		\draw (-100:2) node[fill, circle, inner sep=1.75pt,] {};
		\draw (-140:2) node[fill, circle, inner sep=1.75pt,] {};
	 \draw[ black] (60:2) to[out=60, in=00] (00:2);
	\draw[black] (180:2) to[out=180, in=120] (120:2);
	\draw[black] (-7:2) to[out=-7, in=-60] (-60:2);
	\draw[black,style=dashed] (-100:2) to[out=-100, in=-60] (-60:2);
	\draw[black] (-100:2) to[out=-100, in=-140] (-140:2);
		\draw[black,style=dashed] (-140:2) to[out=-140, in=-180] (-180:2);
	\end{tikzpicture}
\label{fig:path0}
		\end{figure}\\
In the following theorem, we reproduce the proof context outlined in the original proof, employing the same methodology to generate a similar list of sub-digraphs but with a new set of chromatic numbers, enabling us to enhance the bound from $3n$ to $2n$.
		
	\begin{theorem}\label{th1}
	Every $2n$-chromatic Hamiltonian digraph contains a subdivision of every oriented cycle of order $n$.
\end{theorem}
	
\begin{proof}
		Let $D$ be a $2n$-chromatic Hamiltonian digraph and $C$ be a non-directed oriented cycle of order $n$. We aim to show that $D$ contains a subdivision of $C$.\ Consider $k_1,\dots, k_{2t}$ as non-zero positive integers such that $C=C^+(k_1,\dots,k_{2t})$. Given that $k_1 + \dots + k_{2t} = n$, we must have either $k_1 + k_3 + \dots + k_{2t-1} \leq \frac{n}{2}$ or $k_2 + k_4 + \dots + k_{2t} \leq \frac{n}{2}$. Notably, $C^+(k_1, \dots, k_{2t})=C^+(k_{2t},k_{2t-1},\dots,k_1)$. Thus, after relabeling if necessary, we may assume that $k_2 + k_4 + \dots + k_{2t} \leq \frac{n}{2}$. Concerning the number of blocks, the highest possible value for $2t$ is $n$, which occurs when the cycle is antidirected.\ Consequently, $\chi(D)=2n=n+\frac{n}{2}+\frac{n}{2}
		\geq (k_1+\dots +k_{2t})+(k_2+k_4+\dots +k_{2t})+t
		= (k_1-1)+(2k_2+2)+\dots+ (k_{2t-1}-1)+(2k_{2t}+2)$.
		
		Consider a Hamiltonian circuit $C'$ of $D$,  and proceed to construct a sequence of induced sub-digraphs $(D_i)_{1\leq i\leq 2t}$ based on the methodology depicted in Figure \ref{fig:path0} and adapted from the approach introduced in the original proof  by El Joubbeh \cite{El.joubbeh.cycles} (See also the detailed construction in the proof of Theorem \ref{MAmain}). In our adapted construction, induced sub-digraphs are created with the distinction that $\chi(D_i)=k_i-1$ when $i$ is odd and $2k_i+2$ when $i$ is even. In cases where $i$ is odd and $k_i=1$, we treat $D_i$ as an empty digraph. In such instances, after generating $D_{i-1}$, the generation of $D_{i+1}$ begins.
			
	When $i$ is even, according to Corollary \ref{coro1}, $D_{i}$ contains two arcs, $e_{i}$ and $f_{i}$, as $k_{i}$-secant arcs with respect to the natural enumeration induced by $C'[V(D_{i})]$. Without loss of generality, we assume that the arc $e_{i}$ precedes the arc $f_{i}$ in the direction of $C'[V(D_{i})]$. For each $a \in \{e_{i}, f_{i}\}$, we denote the two vertices of $a$ as $x(a)$ and $y(a)$ so that $C'_{[x(a),y(a)]} \subseteq C'[V(D_{i})]$. On the other hand, observe that whenever $i$ is odd, $D_i$ has at least $k_i-1$ vertices, and it is preceded by a vertex and succeeded by a vertex on the cycle $C'$, making a block of length at least $k_i$. Thus, $e_2 \cup C'{[x(f_2),y(e_2)]} \cup f_2 \cup C'{[y(f_2),x(e_4)]} \cup e_4 \cup \dots \cup e_{2t} \cup C'{[x(f_{2t}),y(e_{2t})]} \cup f_{2t} \cup C'{[y(f_{2t}),x(e_2)]}$ is a subdivision of $C^+(k_1,\dots,k_{2t})$ (See Figure \ref{fig:path01}).
		\begin{figure}[h]
		\centering
		\caption{}
		\begin{tikzpicture}
		
			 \draw[line width=1.5pt, red] (-5:2) arc (-5:-37:2);
			  \draw[line width=0.75pt, dashed, brown] (-5:2) arc (-5:15:2);
			  \draw[line width=0.75pt, dashed, brown] (315:2) arc (315:295:2);
			  \draw[line width=0.75pt, dashed, brown] (190:2) arc (190:170:2);
			  \draw[line width=0.75pt, dashed, brown] (135:2) arc (135:110:2);
			 \draw[line width=1.5pt, red] (170:2) arc (170:143:2);
			 \draw[line width=1.5pt, blue] (110:2) arc (110:23:2);
			 \draw[line width=1.5pt, blue] (263:2) arc (263:295:2);
			 \draw[line width=1.5pt, blue] (198:2) arc (198:220:2);
			\node[circle, fill=black, inner sep=1.75pt] (v1) at (15:2) {};
			\node[circle, fill=blue, inner sep=0.6pt] (v21) at (24:2) {};
			\Edge[Direct, color=blue](v21)(v1)
		\node[circle, fill=black, inner sep=1.75pt] (v2) at (-45:2) {};
		\node[circle, fill=red, inner sep=0.6pt] (v22) at (-36:2) {};
		\Edge[Direct, color=red](v22)(v2)
		
		\node[circle, fill=black, inner sep=1.75pt] (v3) at (-5:2) {};
		\node[circle, fill=black, inner sep=1.75pt] (v4) at (-65:2) {};
		\Edge[Direct,bend=80, color=blue](v1)(v2)
			\Edge[Direct,bend=-50,color=red](v2)(v1)
			\Edge[Direct,bend=-60, color=blue](v3)(v4)
			\Edge[Direct,bend=30,color=red](v4)(v3)
			
			\node[circle, fill=black, inner sep=1.75pt] (v5) at (190:2) {};
			\node[circle, fill=blue, inner sep=0.5pt] (v51) at (199:2) {};
			\Edge[Direct, color=blue](v51)(v5)
			\node[circle, fill=black, inner sep=1.75pt] (v6) at (135:2) {};
			\node[circle, fill=red, inner sep=0.5pt] (v61) at (144:2) {};
			\Edge[Direct, color=red](v61)(v6)
			
			\node[circle, fill=black, inner sep=1.75pt] (v7) at (170:2) {};
			\node[circle, fill=black, inner sep=1.75pt] (v8) at (110:2) {};
			\Edge[Direct,bend=80, color=blue](v5)(v6)
			\Edge[Direct,bend=-50,color=red](v6)(v5)
			\Edge[Direct,bend=-60, color=blue](v7)(v8)
			\Edge[Direct,bend=30,color=red](v8)(v7)
			\node[circle, fill=black, inner sep=1.75pt] (v9) at (220:2) {};
			\node[circle, fill=black, inner sep=1.75pt] (v10) at (255:2) {};
			\node[circle, fill=blue, inner sep=0.6pt] (v101) at (264:2) {};
			\Edge[Direct, color=blue](v101)(v10)
				\Edge[Direct, style=dashed,bend=30,color=red](v9)(v10)
				\Edge[Direct, style=dashed,bend=50,color=blue](v10)(v9)
					\node[circle, fill=black, inner sep=1.6pt] (v11) at (237.5:2) {};
					\node[circle, fill=black, inner sep=1.6pt] (v12) at (231:2) {};
					\node[circle, fill=black, inner sep=1.6pt] (v13) at (244:2) {};
					\node[align=center] at (0.8,-0.5) {\footnotesize{$f_2$}};
					\node[align=center] at (2.5,-0.9) {\footnotesize{$e_2$}};
					\node[align=center] at (0.9,-2.05) {\tiny{$y(f_2)$}};
						\node[align=center] at (1.4,-1.65) {\tiny{$y(e_2)$}};
						\node[align=center] at (1.7,0) {\tiny{$x(f_2)$}};
						\node[align=center] at (2.3,0.65) {\tiny{$x(e_2)$}};
						\node[align=center] at (-1.5,-0.35) {\tiny{$x(e_{2t})$}};
						\node[align=center] at (-1.65,0.15) {\tiny{$x(f_{2t})$}};
						\node[align=center] at (-1.15,1.6) {\tiny{$y(e_{2t})$}};
					\node[align=center] at (-0.8,0.5) {\footnotesize{$f_{2t}$}};
					\node[align=center] at (-2.5,0.9) {\footnotesize{$e_{2t}$}};
					\node[align=center] at (-0.7,2.1) {\tiny{$y(f_{2t})$}};
			\end{tikzpicture}
		\label{fig:path01}
	\end{figure}\\
It can be checked that the construction of the desired cycle is valid, regardless of the orientation of any arc $a\in \{e_2,f_2,\dots,e_{2t},f_{2t}\}$.
\end{proof}

\section{Subdivision of oriented cycles in digraphs with Hamiltonian path}\label{H.P.}
The subsequent lemma, well-known in the field of graph theory and used as a key element in proving our main result, imposes an upper bound on the chromatic number of a graph with two spanning subgraphs that together include all its edges.
\begin{lemma} \label{H1H2}
	Consider a graph $G$ with spanning subgraphs $G_1$ and $G_2$ such that $E(G) = E(G_1) \cup E(G_2)$. Then, $\chi(G)\leq \chi(G_1)\chi(G_2)$.
\end{lemma}

\begin{proof}
	For $i \in \{1, 2\}$ let $c_i$ be a proper $\chi(G_i)$-coloring of $G_i$. Now, the function $\Gamma: V(G) \to \{1, \dots, \chi(G_1)\} \times \{1, \dots, \chi(G_2)\}$, defined by $\Gamma(x) = (c_1(x), c_2(x))$ is a proper coloring of $G$. To see this, consider $xy \in E(G)$, then either $xy \in E(G_1)$ or $xy \in E(G_2)$. Thus, $c_1(x) \neq c_1(y)$ or $c_2(x) \neq c_2(y)$, leading to $\Gamma(x) \neq \Gamma(y)$. Hence, $\chi(G) \leq \chi(G_1)\chi(G_2)$.
\end{proof}
This lemma also applies to any digraph containing two subdigraphs that collectively cover all its arcs.
	\begin{theorem}\label{MAmain}
		Every ($12n-5$)-chromatic digraph with a Hamiltonian directed path  contains a subdivision of every oriented cycle of order $n$ with more than one block.
	\end{theorem}
\begin{proof}
		Let $H$ be a ($12n-5$)-chromatic digraph with a Hamiltonian directed path and let $P$ be a directed path of $H$. Let $C$ be a non-directed oriented cycle of order $n$, and assume that $H$ contains no subdivision of $C$.\ Form a list of induced subdigraphs of $H$, $(H_i)_{1\leq i\leq \ell}$ in the following way. Starting  from the first vertex of $P$ and adding vertices one successor at a time until we reach a vertex where the induced subdigraph of $H$ formed by all these vertices achieves a chromatic number of $2n$.\ Omit the last added vertex and label the induced subdigraph formed by the preceding vertices as $H_1$.\ This process is feasible because whenever a new vertex is added, the chromatic number either remains the same or increases by one.\ Thus, $\chi(H_1)=2n-1$, and if we add the next successor the chromatic number will jump to $2n$. Continue from the successor of the last vertex in $P[V(H_1)]$ to construct $D_2$ using the same method as employed for $H_1$. Repeat this process to generate $H_3, H_4,$ and so on until all vertices of $H$ are exhausted, with $\ell$ representing the index of the last created subdigraph (See Figure \ref{fig:path}).
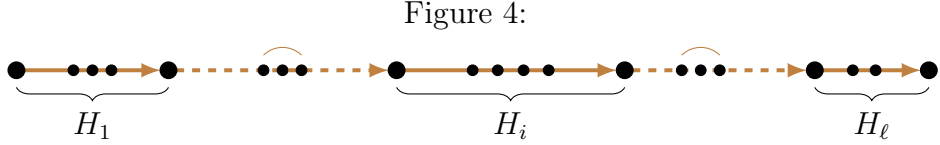
	\begin{figure}[h] 
		\centering
\caption{}
	
	\begin{tikzpicture}
		\Vertex[size=0.2,color=black,position=below]{v1}
		\Vertex[x=2,size=0.2,color=black,position=below]{v2}
		\Vertex[x=5,size=0.2,color=black,position=below]{v3}
		\Vertex[x=8,size=0.2,color=black,position=below]{v4}
			\Vertex[x=10.5,size=0.2,color=black,position=below]{v5}
		\Vertex[x=12,size=0.2,color=black,position=below]{v6}
		\Edge[Direct,color=brown](v1)(v2)
		\Edge[Direct,color=brown](v3)(v4)
		\Edge[Direct,color=brown](v5)(v6)
		\Edge[Direct,style=dashed,color=brown](v2)(v3)
		\Edge[Direct,style=dashed,color=brown](v4)(v5)
	
		\draw[color=brown] (3.25,0.2) to[out=45, in=135] (3.75,0.2);
		\draw[color=brown] (8.75,0.2) to[out=45, in=135] (9.25,0.2);
		\foreach \x in {0.75,1,1.25,3.25,3.5,3.75,6,6.33,6.67,7,8.75,9,9.25,11,11.3}{
			\filldraw[black] (\x,0) circle (2pt);
		}
	\draw [decorate,decoration={brace,amplitude=5pt,mirror,raise=4ex}]
	(0,0.5) -- (2,0.5) node[midway,yshift=-3em]{$H_1$};
	\draw [decorate,decoration={brace,amplitude=5pt,mirror,raise=4ex}]
	(5,0.5) -- (8,0.5) node[midway,yshift=-3em]{$H_i$};
	\draw [decorate,decoration={brace,amplitude=5pt,mirror,raise=4ex}]
	(10.5,0.5) -- (12,0.5) node[midway,yshift=-3em]{$H_\ell$};
	\end{tikzpicture}

\label{fig:path}
\end{figure}\\
 Note that the chromatic number of $H_\ell$ need not reach $2n-1$. Also, the property seen in $H_1$ stays the same for all $H_i$ up to $\ell-1$. That is, $\chi(H_i)=2n-1$, and adding the first vertex from $P[H_{i+1}]$ causes the chromatic number to instantly reach $2n$.\\
 
Set $R_1:=(V(H);\bigcup_{i=1}^\ell E(H_i))$ and $R_2:=(V(H);E(H)\setminus\bigcup_{i=1}^\ell E(H_i))$. Then these two spanning subdigraphs cover all the arcs of $H$. Thus, by Lemma \ref{H1H2}, we have: $$\chi(H)\leq \chi(R_1)\chi(R_2).\ \ (\star)$$

 The digraph $R_1$ comprises $\ell$ digraphs, $H_1,\dots,H_{\ell}$. Each of these has a chromatic number of $2n-1$, except the last one, which possesses a chromatic number less than or equal to $2n-1$. These digraphs are vertex-disjoint, and within $R_1$, there are no arcs linking any two of them. Consequently, $\chi(R_1)=\max\{\chi(H_i) \ | \ 1\leq i\leq \ell\}=2n-1$.\\
  
  For $R_2$, we claim that $\chi (R_2)\leq 6$, which would finish the proof by establishing a contradiction when using Display ($\star$).
  
    \begin{claim}\label{claim2}
  	$\chi(R_2)\leq 6$.
  \end{claim}
 \textit{Proof of Claim \ref{claim2}.} Consider the graph $G$ obtained from $R_2$ by contracting all vertices in $V(H_i)$, $1\leq i\leq \ell$, into a single vertex labeled $v_i$, where $v_rv_s\in E(G)$ if and only if there exist $x\in V(H_r)$ and $y\in V(H_s)$ such that $(x,y)\in E(R_2)$ or $(y,x)\in E(R_2)$. Based on the construction of $R_2$, $V(H_i)$ is a stable set in $R_2$, $ \forall \ 1\leq i \leq \ell $. Consequently, any proper coloring $c'$ of $G$ can generate a proper coloring $c$ of $R_2$ by assigning all vertices in $V(H_i)$ in $c$ the color $c'(v_i)$. Therefore, $\chi(R_2)\leq \chi(G)$.\\
 
  Now, set $S_1=\{v_i \ | \ i \ \textnormal{is odd}\}$, $S_2=\{v_i \ | \ i \ \textnormal{is even}\}$, $G_1=G[S_1]$ and $G_2=G[S_2]$. Thus, $\chi(G)\leq \chi(G_1)+\chi(G_2)$. Set $N_1=(v_1,v_3,\dots,v_{\ell_1})$ and $N_2=(v_2,v_4,\dots,v_{\ell_2})$, where $\ell_1$ and $\ell_2$ represent the largest odd and even integers less than or equal to $\ell$, respectively. Thus, $N_1$ and $N_2$ are two enumerations of $G_1$ and $G_2$, respectively.\\
  
   We argue that $\chi(G_1)\leq 3$. Assuming that $\chi(G_1)\geq 4$, then by Lemma \ref{L:secant} applied with $k=1$, $G_1$ contains $1$-secant edges $v_iv_j$ and $v_rv_s$, $i<r< j< s$, with respect to $N_1$. Thus, there exist $u_i\in V(H_i)$, $u_j\in V(H_j)$, $u_r\in V(H_r)$ and $u_s\in V(H_s)$ such that ``$(u_i, u_j)$ or $(u_j,u_i)$'' and ``$(u_r, u_s)$ or $(u_s, u_r)$'' are arcs in $R_2$. Note that $v_i$, $v_j$, $v_r$ and $v_s$ are  vertices within $S_1$, implying that $i$, $j$, $r$, and $s$ are all odd integers, it follows that between them, there exist two even integers, say $m_1$ and $m_2$, such that $i< m_1< r $ and $j < m_2< s$. We study two cases:\\

\noindent\textbf{Case 1.} If $(u_j,u_i)\in E(R_2)$ or $(u_s,u_r)\in E(R_2)$. These two scenarios are analogous. Assuming that $(u_j,u_i)\in E(H_2)$. Set $D=H[V(P_{[u_i,u_j]})]$. Thus $D$ is a Hamiltonian digraph. Additionally, $V(H_{m_1})\cup\{x\}\subseteq V(D)$, where $x$ represents the starting vertex of $P[V(H_{m_1+1})]$ (See Figure \ref{fig2:path}).
	\begin{figure}[h] 
	\centering
	\caption{}
	
	\begin{tikzpicture}
		\Vertex[x=8.3,size=0.2,color=black,label=$x$,position=below]{v14}
		\Vertex[size=0.2,color=black,position=below]{v1}
		\Vertex[x=2,size=0.2,color=black,position=below]{v2}
		\Vertex[x=5,size=0.2,color=black,position=below]{v3}
		\Vertex[x=8,size=0.2,color=black,position=below]{v4}
		\Vertex[x=10.5,size=0.2,color=black,position=below]{v5}
		\Vertex[x=12,size=0.2,color=black,position=below]{v6}
		\Vertex[x=4,size=0.2,color=black,label=$u_i$,position=below]{v7}
		\Vertex[x=9,size=0.2,color=black,label=$u_j$,position=below]{v8}
		\Edge[Direct,color=brown](v1)(v2)
		\Edge[Direct,color=brown](v3)(v4)
		\Edge[Direct,color=brown](v5)(v6)
		\Edge[Direct,style=dashed,color=brown](v2)(v3)
		\Edge[Direct,style=dashed,color=brown](v4)(v5)
		\Edge[Direct,bend=-40, color=black](v8)(v7)
		
		\draw[color=brown] (2.75,0.2) to[out=45, in=135] (3.25,0.2);
		\draw[color=brown] (9.25,0.2) to[out=45, in=135] (9.75,0.2);
		\foreach \x in {0.75,1,1.25,2.75,3,3.25,6,6.33,6.67,7,9.25,9.5,9.75,11,11.3}{
			\filldraw[black] (\x,0) circle (2pt);
		}
		\draw [decorate,decoration={brace,amplitude=5pt,raise=-1.8ex}]
		(5,0.5) -- (8,0.5) node[midway,yshift=0.2em]{$H_{m_1}$};
		\draw [decorate,decoration={brace,amplitude=6pt,mirror,raise=6ex}]
		(4,0.6) -- (9,0.6) node[midway,yshift=-3.8em]{$D$};
	\end{tikzpicture}
	
	\label{fig2:path}
\end{figure}\\
 Recall that the properties of $H_{i}$, $1\leq i\leq \ell-1$, not only relate to a chromatic number of $2n-1$ but also, adding the first vertex of $P[V(H_{i+1})]$ increases the chromatic number to $2n$.  This yields $\chi (D)\geq \chi(H[V(H_{m_1})\cup\{x\}])= 2n$. According to Theorem \ref{th1}, this implies the existence of a subdivision of $C$ within $D$, and consequently within $H$, a contradiction.\\

\noindent\textbf{Case 2.} If $(u_i,u_j)\in E(R_2)$ and $(u_r,u_s)\in E(R_2)$. Set $C=C^+(k_1,\dots,k_{2t})$. Following the proof of Theorem \ref{th1}, it is reasonable to presume that $k_2+k_4+\dots +k_{2t}\leq \frac{n}{2}$, and $\chi(H_{m_1})=2n-1\geq [(k_1-1)+(2k_2+2)+\dots+ (k_{2t-1}-1)+(2k_{2t}+2)]-1$.\ Utilizing the ``vertex addition method'' previously employed in constructing $(H_i)_{1\leq i \leq \ell}$, the chromatic number of $H_{m_1}$ is adequate to generate a list of $2t-1$ subdigraphs ``$(D_i)_{1\leq i \leq 2t-1}$'' defined as follows:  Start from the initial vertex within $P[V(H_{m_1})]$, add its successor vertex on $P$, and iterate until achieving an induced subdigraph of $H_{m_1}$ with a chromatic number of $k_1-1$, then label it $D_1$.\ Progress from the subsequent vertex after $P[V(D_1)]$ to form $D_2$, stopping upon reaching a chromatic number of $2k_2+2$. \ Continue this process, setting $\chi(D_i)=k_i-1$ when $i$ is odd, and $\chi(D_i)=2k_i+2$ when $i$ is even. Repeat until $D_{2t-1}$ is obtained, merging the remaining vertices of $H_{m_1}$ into $D_{2t-1}$ (See Figure \ref{fig4:path}).
\begin{figure}[h] 
	\centering
	\caption{}
	
	\begin{tikzpicture}
		\Vertex[size=0.2,color=black,position=below]{v1}
		\Vertex[x=2,size=0.2,color=black,label=$u_i$,position=above]{v2}
		\Vertex[x=2.95,size=0.2,color=black,,position=below]{v11}
		\Vertex[x=4.4,size=0.2,color=black,,position=below]{v12}
		\Vertex[x=6,size=0.2,color=black,,position=below]{v13}
		\Vertex[x=7.25,size=0.2,color=black,,position=below]{v14}
		\Vertex[x=8,size=0.2,color=black,label=$u_r$,position=below]{v3}
		\Vertex[x=9,size=0.2,color=black,label=$u_j$,position=below]{v4}
		\Vertex[x=10,size=0.2,color=black,label=$u_s$,position=below]{v5}
		\Vertex[x=12,size=0.2,color=black,position=below]{v6}
		\Edge[Direct,style=dashed,color=brown](v1)(v2)
		\Edge[Direct,style=dashed,color=brown](v2)(v11) 
		\Edge[Direct,,style=dashed, color=brown](v3)(v4)
		\Edge[Direct,style=dashed, color=brown](v5)(v6)
		\Edge[Direct, color=brown](v11)(v12) 
		\Edge[Direct,color=brown](v13)(v14)
		\Edge[Direct, style=dashed,color=brown](v12)(v13)
		\Edge[Direct,style=dashed, color=brown](v14)(v3)
		\Edge[Direct,style=dashed,color=brown](v4)(v5)
		\Edge[Direct,bend=35, color=black](v2)(v4)
		\Edge[Direct,bend=-50, color=black](v3)(v5)
		
			\draw [decorate,decoration={brace,amplitude=5pt,raise=-1.8ex}]
		(2.95,0.5) -- (7.25,0.5) node[midway,yshift=0.2em]{$H_{m_1}$};
		
		\draw[color=brown] (0.6,0.2) to[out=45, in=135] (1.2,0.2);
		\draw[color=brown] (10.6,0.2) to[out=45, in=135] (11.2,0.2);
		\foreach \x in {0.6,0.9,1.2,3.3,3.6,3.9,5.15, 4.95, 5.35,6.4,6.7,10.6,10.9,11.2}{
			\filldraw[black] (\x,0) circle (2pt);
		}
		\draw [decorate,decoration={brace,amplitude=5pt,mirror,raise=4ex}]
		(2.95,0.5) -- (4.4,0.5) node[midway,yshift=-3em]{$D_1$};
			\draw [decorate,decoration={brace,amplitude=5pt,mirror,raise=4ex}]
		(6,0.5) -- (7.25,0.5) node[midway,yshift=-3em]{$D_{2t-1}$};
	\end{tikzpicture}
	
	\label{fig4:path}
\end{figure}\\
Following the completion of the above list generation, we could use Corollary \ref{coro1} to observe that whenever $i$ is even, $D_{i}$ contains two arcs $e_{i}$ and $f_{i}$ as $k_{i}$-secant arcs with respect to the natural enumeration induced by $P[V(D_{i})]$.\ Without loss of generality, suppose that the arc $e_{i}$ comes before the arc $f_{i}$ with respect to the direction of $P[V(D_{i})]$, and for each $a\in \{e_{i},f_{i}\}$, we denote the two vertices of $a$ as $x(a)$ and $y(a)$ in such a way that $x(a)$ precedes $y(a)$ on the directed path $P[V(D_{i})]$.\ On the other hand, observe that whenever $i$ is odd, $D_i$ has at least $k_i-1$ vertices, and it is preceded by a vertex and succeeded by a vertex on the path $P$, making a block of length at least $k_i$.\ Set $P_1=P_{[u_i,x(e_2)]}\cup e_2\cup P_{[x(f_2),y(e_2)]}\cup f_2\cup P_{[y(f_2),x(e_4)]}\cup e_4\cup \dots \cup e_{2t-2}\cup P_{[x(f_{2t-2}),y(e_{2t-2})]}\cup f_{2t-2}\cup P_{[y(f_{2t-2}),u_r]}\cup (u_r,u_s)$ and $P_2=\{u_i\}\cup (u_i,u_j)\cup P_{[u_j,u_s]}$. Consequently, $P_1$ is an oriented $u_i,u_s$-path that forms a subdivision of $P^+(k_1,\dots,k_{2t-1})$, and regarding $P_2$, it is a directed $u_i,u_s$-path, and as $j< m_2< s$, we have $V(H_{m_2})\cup\{u_i,u_j,u_s\}\subseteq V(P_2)$, implying $l(P_2)=|V(P_2)|-1\geq |V(H_{m_2})\cup\{u_i,u_j,u_s\}|-1\geq \chi(H_{m_2})+2 =2n+1 \gg k_{2t}$. Moreover, with $V(P_1)\cap V(P_2)=\{u_i,u_s\}$, the union $P_1\cup P_2=C^+(k'_1,\dots,k'_{2t})$ is a cycle formed, where $k'_i\geq k_i$ holds for all $1\leq i\leq 2t$. (See Figure \ref{fig3:path}). Then $P_1\cup P_2$ it is a subdivision of $C$, a contradiction.\\
\begin{figure}[h] 
	\centering
	\caption{}
	
	\begin{tikzpicture}
		\Vertex[size=0.2,color=black,position=below]{v1}
		\Vertex[x=1,size=0.25,color=black,label=$u_i$,position=below]{v2}
			\Vertex[x=2.2,size=0.2,color=black]{v9}
			\Vertex[x=3.9,size=0.2,color=black,]{v10}
			\Vertex[x=2.95,size=0.2,color=black]{v11}
			\Vertex[x=4.7,size=0.2,color=black]{v12}
			\Vertex[x=5.65,size=0.2,color=black,,position=below]{v13}
			\Vertex[x=6.75,size=0.2,color=black,,position=below]{v14}
		\Vertex[x=8,size=0.2,color=black,label=$u_r$,position=below]{v3}
		\Vertex[x=9,size=0.2,color=black,label=$u_j$,position=below]{v4}
		\Vertex[x=11,size=0.25,color=black,label=$u_s$,position=below]{v5}
		\Vertex[x=12,size=0.2,color=black,position=below]{v6}
		\Edge[Direct,style=dashed,color=brown](v1)(v2)
		\Edge[Direct,style=dashed,color=brown](v9)(v11) 
		\Edge[Direct,style=dashed,color=brown](v10)(v12) 
		\Edge[Direct,,style=dashed, color=brown](v3)(v4)
		\Edge[Direct,style=dashed, color=brown](v5)(v6)
		\Edge[Direct, color=red](v11)(v10)
		\Edge[Direct, color=blue,bend=60](v11)(v12) 
		\Edge[Direct, color=red,bend=-100](v12)(v11) 
		\Edge[Direct, color=blue](v12)(v13)
		\Edge[Direct, style=dashed, bend=35, color=blue](v13)(v14)
		\Edge[Direct, style=dashed, bend=35,  color=red](v14)(v13)
	\node[align=center] at (2.3,0.25) {\tiny{$ x(e_2)$}};
	\node[align=center] at (3.05,-0.25) {\tiny{$ x(f_2)$}};
	\node[align=center] at (4,0.25) {\tiny{$ y(e_2)$}};
	\node[align=center] at (4.8,-0.25) {\tiny{$ y(f_2)$}};
	\node[align=center] at (3.05,-0.8) {\footnotesize{$e_2$}};
	\node[align=center] at (4,0.8) {\footnotesize{$f_2$}};
	\node[align=center] at (5.35,-0.5) {$P_1$};
	\node[align=center] at (8.75,0.5) {$P_2$};
		\Edge[Direct, color=blue](v14)(v3)
		\Edge[Direct,color=blue](v2)(v9)
		\node[align=center] at (3.4,0.2) {\scriptsize{\color{red}$k'_2$}};
		\node[align=center] at (7.3,0.2) {\scriptsize{\color{blue}$k'_{2t-1}$}};
		\node[align=center] at (10,0.2) {\scriptsize{\color{red}$k'_{2t}$}};
		\node[align=center] at (5.1,0.2) {\scriptsize{\color{blue}$k'_3$}};
		\node[align=center] at (1.7,0.2) {\scriptsize{\color{blue}$k'_1$}};
		\Edge[Direct,bend=95,color=red, color=red](v10)(v9)
		\Edge[Direct,bend=-60,color=blue](v9)(v10)
		\Edge[Direct,color=red](v4)(v5)
		\Edge[Direct,bend=35, color=red](v2)(v4)
		\Edge[Direct,bend=-45,color=blue](v3)(v5)
		
		\foreach \x in {6, 6.2,6.4}{
			\filldraw[black] (\x,0) circle (2pt);
		}
	\end{tikzpicture}
	
	\label{fig3:path}
\end{figure}\\
 The contradictions in \textbf{Case 1} and \textbf{Case 2} confirm that $\chi(G_1)\leq 3$, and similarly $\chi(G_2)\leq 3$. Thus, $\chi(R_2)\leq \chi(G)\leq \chi(G_1)+\chi(G_2)\leq 6$. This completes the proof of Claim \ref{claim2}.

\end{proof}

	\renewcommand{\abstractname}{Acknowledgements}
	\begin{abstract}
		Research of the first author was partially supported by the University Research Board of the American University of Beirut (Project Number 27190.)
	\end{abstract}

\end{document}